\documentclass{amsart}

\title{How the continuum hypothesis could have been a fundamental axiom}
\usepackage[pdfauthor={Joel David Hamkins},
    pdftitle={How the continuum hypothesis could have been a fundamental axiom},
    hidelinks]{hyperref}

\author{Joel David Hamkins}
\address[Joel David Hamkins]
{O'Hara Professor of Logic, University of Notre Dame, 100 Malloy Hall, Notre Dame, IN 46556 USA, \&\ V. Research Fellow, Mathematical Institute, University of Oxford, UK.}
\email{jdhamkins@nd.edu}
\urladdr{http://jdh.hamkins.org}

\thanks{Commentary on this article can be made on the author's blog at \href{https://jdh.hamkins.org/how-ch-could-have-been-fundamental}{https://jdh.hamkins.org/how-ch-could-have-been-fundamental}. This article follows up on a line of argument I had presented in~\cite[p.78-79,~Questions~2.18,~2.19,~2.20]{Hamkins2021:Lectures-on-the-philosophy-of-mathematics}. I am grateful to the anonymous referees for several helpful suggestions, as well as to Dan Isaacson, W.~Hugh Woodin, Mel Fitting, Jason Chen, Tim Button, as well as to other participants at the talks I gave on this topic at the 2024 Eastern APA in New York, at UC Irvine, and at the University of Oslo.}

\usepackage{latexsym,amsfonts,amsmath,amssymb,mathrsfs}
\usepackage[dvipsnames,svgnames]{xcolor}
\usepackage{wrapfig} 
\usepackage{caption} 
\DeclareCaptionStyle{compact}{font=footnotesize,format=plain,name={Fig},labelsep=period,skip=1.5ex,margin=0pt,oneside,justification=centering}
\DeclareCaptionStyle{leftside}{style=compact,margin={0pt,9pt}}
\DeclareCaptionStyle{rightside}{style=compact,margin={9pt,0pt}}
\usepackage{tikz}
\usetikzlibrary{arrows,arrows.meta,petri,topaths,positioning,shapes,shapes.misc,patterns,calc,decorations.pathreplacing,hobby}
\usepackage[shortlabels]{enumitem} 

\newtheorem*{theorem*}{Theorem}

\newtheorem*{maintheorem*}{Main Theorem}
\newtheorem*{maintheorems*}{Main Theorems}

\newtheorem*{corollary*}{Corollary}
\newtheorem*{corollaries*}{Corollaries}

\theoremstyle{definition}

\newtheorem*{definition*}{Definition}

\newtheorem*{question*}{Question}

\newtheorem*{questions*}{Questions}
\newtheorem*{mainquestion*}{Main Question} 
\newtheorem*{openquestion*}{Open Question} 
\theoremstyle{remark}

\theoremstyle{theorem}
\newcommand{\QED}{\end{proof}}

\def\proclaim[#1]{{\bf #1}}
\def\BF#1.{{\bf #1.}}

\def\says#1:#2\par{\item[#1] #2\par}

\newcommand{\Godel}{G\"odel}

\newcommand{\C}{{\mathbb C}}

\newcommand{\N}{{\mathbb N}}

\newcommand{\Q}{{\mathbb Q}}

\newcommand{\Z}{{\mathbb Z}}
\newcommand{\R}{{\mathbb R}}

\newcommand{\continuum}{\mathfrak{c}}

\makeatletter
\newcommand{\dotminus}{\mathbin{\text{\@dotminus}}}
\newcommand{\@dotminus}{%
  \ooalign{\hidewidth\raise1ex\hbox{.}\hidewidth\cr$\m@th-$\cr}%
}
\makeatother

\newcommand{\of}{\subseteq}

\newcommand{\elesub}{\prec}

\DeclareMathOperator{\necessary}{\text{\tikz[scale=.6ex/1cm,baseline=-.6ex,line width=.1ex]{\draw (-1,-1) rectangle (1,1);}}}

\newcommand\dbrace{\hskip-1.5em\raise3pt\hbox{\rotatebox[origin=c]{-35}{$\left.\strut^{\phantom{|}}\right\}$}}}

\newcommand\UParroW{{\setbox0\hbox{$\Uparrow$}\rlap{\hbox to \wd0{\hss$\mid$\hss}}\box0}}

\renewcommand{\setminus}{\raise.3ex\hbox{\rotatebox{-20}{$-$}}} 

\newcommand{\smalllt}{\mathrel{\mathchoice{\raise2pt\hbox{$\scriptstyle<$}}{\raise1pt\hbox{$\scriptstyle<$}}{\raise0pt\hbox{$\scriptscriptstyle<$}}{\scriptscriptstyle<}}}
\newcommand{\smallleq}{\mathrel{\mathchoice{\raise2pt\hbox{$\scriptstyle\leq$}}{\raise1pt\hbox{$\scriptstyle\leq$}}{\raise1pt\hbox{$\scriptscriptstyle\leq$}}{\scriptscriptstyle\leq}}}

   \def\DHLhksqrt#1#2{%
   \setbox0=\hbox{$#1\sqrt{#2\,}$}\dimen0=\ht0
   \advance\dimen0-0.2\ht0
   \setbox2=\hbox{\vrule height\ht0 depth -\dimen0}%
   {\box0\lower0.4pt\box2}}

\def\[#1]{\mathopen{\lbrack\!\lbrack}#1\mathclose{\rbrack\!\rbrack}}
\newbox\gnBoxA
\newbox\gnBoxB
\newdimen\gnCornerHgt
\setbox\gnBoxA=\hbox{\tiny$\ulcorner$}
\global\gnCornerHgt=\ht\gnBoxA
\newdimen\gnArgHgt
\def\gcode #1{%
\setbox\gnBoxA=\hbox{$#1$}%
\setbox\gnBoxB=\hbox{$\bar #1$}%
\gnArgHgt=\ht\gnBoxB%
\ifnum     \gnArgHgt<\gnCornerHgt \gnArgHgt=0pt%
\else \advance \gnArgHgt by -\gnCornerHgt%
\fi \raise\gnArgHgt\hbox{\tiny$\ulcorner$} \box\gnBoxA %
\raise\gnArgHgt\hbox{\tiny$\urcorner$}}
\newcommand{\UnderTilde}[1]{{\setbox1=\hbox{$#1$}\baselineskip=0pt\vtop{\hbox{$#1$}\hbox to\wd1{\hfil$\sim$\hfil}}}{}}
\newcommand{\Undertilde}[1]{{\setbox1=\hbox{$#1$}\baselineskip=0pt\vtop{\hbox{$#1$}\hbox to\wd1{\hfil$\scriptstyle\sim$\hfil}}}{}}
\newcommand{\undertilde}[1]{{\setbox1=\hbox{$#1$}\baselineskip=0pt\vtop{\hbox{$#1$}\hbox to\wd1{\hfil$\scriptscriptstyle\sim$\hfil}}}{}}
\newcommand{\UnderdTilde}[1]{{\setbox1=\hbox{$#1$}\baselineskip=0pt\vtop{\hbox{$#1$}\hbox to\wd1{\hfil$\approx$\hfil}}}{}}
\newcommand{\Underdtilde}[1]{{\setbox1=\hbox{$#1$}\baselineskip=0pt\vtop{\hbox{$#1$}\hbox to\wd1{\hfil\scriptsize$\approx$\hfil}}}{}}

\renewcommand{\iff}{\mathrel{\leftrightarrow}}

\def\<#1>{\left\langle#1\right\rangle}

\newcommand\No{\N{\text{o}}}

\newcommand{\ZFC}{{\rm ZFC}}
\newcommand{\ZF}{{\rm ZF}}

\newcommand{\CH}{{\rm CH}}

\newcommand{\GCH}{{\rm GCH}}

\newcommand{\MM}{{\rm MM}}

\newcommand{\PFA}{{\rm PFA}}

\newcommand{\MP}{{\rm MP}}

\newcommand{\cell}[1]{\boxit{\hbox to 17pt{\strut\hfil$#1$\hfil}}}
\newcommand{\head}[2]{\lower2pt\vbox{\hbox{\strut\footnotesize\it\hskip3pt#2}\boxit{\cell#1}}}
\newcommand{\boxit}[1]{\setbox4=\hbox{\kern2pt#1\kern2pt}\hbox{\vrule\vbox{\hrule\kern2pt\box4\kern2pt\hrule}\vrule}}
\newcommand{\Col}[3]{\hbox{\vbox{\baselineskip=0pt\parskip=0pt\cell#1\cell#2\cell#3}}}
\newcommand{\tapenames}{\raise 5pt\vbox to .7in{\hbox to .8in{\it\hfill input: \strut}\vfill\hbox to
.8in{\it\hfill scratch: \strut}\vfill\hbox to .8in{\it\hfill output: \strut}}}
\newcommand{\Head}[4]{\lower2pt\vbox{\hbox to25pt{\strut\footnotesize\it\hfill#4\hfill}\boxit{\Col#1#2#3}}}
\newcommand{\Dots}{\raise 5pt\vbox to .7in{\hbox{\ $\cdots$\strut}\vfill\hbox{\ $\cdots$\strut}\vfill\hbox{\
$\cdots$\strut}}}
\usepackage[backend=bibtex,style=alphabetic,maxbibnames=15,maxcitenames=6,dateabbrev=false]{biblatex}
\renewcommand{\UrlFont}{} 
\renewbibmacro{in:}{\ifentrytype{article}{}{\printtext{\bibstring{in}\intitlepunct}}} 
\DeclareFieldFormat{url}{\UrlFont\url{#1}} 
\DeclareFieldFormat{urldate}{
  (version \thefield{urlday}\addspace%
  \mkbibmonth{\thefield{urlmonth}}\addspace%
  \thefield{urlyear}\isdot)}
\DeclareFieldFormat{eprint:arxiv}{
  \ifhyperref
    {\href{http://arxiv.org/abs/#1}{%
        arXiv\addcolon\nolinkurl{#1}}\iffieldundef{eprintclass}{}{\UrlFont{\mkbibbrackets{\thefield{eprintclass}}}}}
    {arXiv\addcolon\nolinkurl{#1}\iffieldundef{eprintclass}{}{\UrlFont{\mkbibbrackets{\thefield{eprintclass}}}}}}
\begin{document}

\begin{abstract}
I describe a simple historical thought experiment showing how we might have come to view the continuum hypothesis as a fundamental axiom, one necessary for mathematics, indispensable even for calculus.
\end{abstract}

\maketitle

\section{Introduction}

I should like to describe how our attitude toward the continuum hypothesis could easily have been very different than it is. If our mathematical history had been just a little different, I claim, if certain mathematical discoveries had been made in a slightly different order, then we would naturally view the continuum hypothesis as a fundamental axiom of set theory, one furthermore necessary for mathematics and indeed indispensable for making sense of the core ideas underlying calculus.

The continuum hypothesis (\CH) is the assertion that the cardinality of the set of real numbers is the first uncountable infinity, or in other words, that $2^{\aleph_0}=\aleph_1$. This hypothesis is known to be independent of the Zermelo-Fraenkel $\ZFC$ axioms of set theory---it is neither provable nor refutable, if $\ZFC$ itself is consistent, and it remains independent even relative to any of the usual large cardinal axioms. Currently~\CH\ is not generally regarded as part of the standard axiomatization of set theory, but rather is taken as a separate supplemental hypothesis, one to be mentioned explicitly, assumed or denied, proved or refuted, in diverse circumstances for different purposes. The continuum hypothesis holds, for example, in the constructible universe $L$ and in other canonical inner models, but we can make it fail (or hold) in forcing extensions; it is refuted by the forcing axioms \PFA\ and \MM, which settle the continuum as $\aleph_2$. In the study of the cardinal characteristics of the continuum, set theorists routinely work with $\neg\CH$ as the subject is trivialized in a sense under $\CH$, since there would be no room for variation in the cardinal characteristics, although it is also trivialized in a different way under the forcing axioms, since these imply that they are all fully pushed up to value continuum.

Since the truth or falsity of~\CH\ cannot be settled on the basis of proof from the $\ZFC$ axioms, set theorists have offered various philosophical arguments aiming at a solution to the continuum problem, the problem of determining whether~\CH\ holds or its negation. (See my survey discussion in~\cite[chapter~8]{Hamkins2021:Lectures-on-the-philosophy-of-mathematics}.) For example, Chris Freiling~\cite{Freiling1986:AxiomsOfSymmetry:ThrowingDartsAtTheRealLine} advances an argument for $\neg\CH$ based on prereflective intuitions about randomness as a primitive notion. W.~Hugh Woodin made a case for $\neg\CH$ based on considerations of $\Omega$-logic and forcing absoluteness (see the survey in~\cite{Koellner2023:SEP-continuum-hypothesis}). More recently, Woodin argues on the other side, making a case for~\CH\ based on features of his theory of Ultimate $L$, a canonical inner model accommodating even the largest large cardinals (see~\cite{Rittberg2015:How-Woodin-changed-his-mind-new-thoughts-on-CH} for an account of Woodin's change of heart). Defending set-theoretic pluralism, I argue in~\cite{Hamkins2012:TheSet-TheoreticalMultiverse} that it is incorrect to describe~\CH\ as an open question---the answer to~\CH, rather, is pluralist, consisting of the deep body of knowledge that we have concerning how it behaves in the set-theoretic multiverse, how we can force it or its negation while preserving diverse other set-theoretic features.

Many set theorists have yearned for what I call the \emph{dream solution} to the continuum problem, by which we settle the~\CH\ once and for all by introducing a new set-theoretic principle, the ``missing'' axiom, which everyone agrees is fully consonant with the concept of set and which also provably settles~\CH. I argue in~\cite{Hamkins2015:IsTheDreamSolutionToTheContinuumHypothesisAttainable}, however, that this will never happen.
\begin{quote}\footnotesize
  Our situation with~\CH\ is not merely that~\CH\ is formally independent and we have no additional knowledge about whether it is true or not. Rather, we have an informed, deep understanding of how it could be that~\CH\ is true and how it could be that~\CH\ fails. We know how to build the~\CH\ and $\neg\CH$ worlds from one another. Set theorists today grew up in these worlds, comparing them and moving from one to another while controlling other subtle features about them. Consequently, if someone were to present a new set-theoretic principle $\Phi$ and prove that it implies $\neg\CH$, say, then we could no longer look upon $\Phi$ as manifestly true for sets. To do so would negate our experience in the~\CH\ worlds, which we found to be perfectly set-theoretic. It would be like someone proposing a principle implying that only Brooklyn really exists, whereas we already know about Manhattan and the other boroughs. And similarly if $\Phi$ were to imply $\CH$. We are simply too familiar with universes exhibiting both sides of~\CH\ for us ever to accept as a natural set-theoretic truth a principle that is false in some of them.
\end{quote}
Nevertheless, I should like to explain in this article how it all might easily have been different. If our mathematical history had been slightly revised in a way I shall presently describe, if certain mathematical discoveries had been made in a slightly different order, then we might have come to look upon~\CH\ as fundamental principle for set theory, one necessary to make sense of mathematical ideas at the core of classical mathematics.

\section{The thought experiment---two number realms}

As a thought experiment, let us imagine that Newton and Leibniz in the early days of calculus provide somewhat fuller accounts of their ideas about infinitesimals. In the actual world, to be sure, a satisfactory account of the basic nature and features of infinitesimals was lacking---the foundations of calculus were famously mocked by Berkeley~\cite{Berkeley1734:A-discourse-addressed-to-an-infidel-mathematician} with withering criticism:
\begin{quote}\footnotesize
  And what are these same evanescent Increments? They are neither finite Quantities nor Quantities infinitely small, nor yet nothing. May we not call them the ghosts of departed quantities?
\end{quote}
It was simply not clear enough in the early accounts of calculus what kind of thing the infinitesimals were and whether they were part of the ordinary number system or somehow transcending it, inhabiting a different larger realm of numbers.\goodbreak

According to Jesseph~\cite[p.168]{Jesseph1993:Berkeleys-philosophy-of-mathematics}, Berkeley argued that
\begin{quote}\footnotesize
If infinitesimal magnitudes are introduced into analysis, the question arises whether they obey the ordinary laws of addition, subtraction, multiplication, and division.%
\end{quote}
And he found fault with both sides of the resulting dichotomy.\goodbreak

What I propose is that we imagine that Newton and Leibniz provide greater clarity concerning the conception of infinitesimals. Specifically, I would like to imagine that Newton and Leibniz conceive of the infinitesimals, as many do today, as living in a larger field of numbers, distinct from but extending the ordinary real numbers. Let us suppose that they
posit two ``realms'' of numbers, the ordinary realm $\R$ of the real numbers and a further realm $\R^*$ consisting of what we might call the \emph{hyperreal} numbers, to use the contemporary terminology, a transcendent number field accommodating the infinitesimals.

This idea alone, that infinitesimals inhabit another realm of numbers, immediately addresses the mocking Berkeley criticism, releasing the tension of the otherwise paradoxical claim that infinitesimals are positive yet also smaller than every positive number, for we need only claim that infinitesimals are smaller than every positive real number, of course, and not smaller than all the other infinitesimal numbers or themselves. The two-number-realms idea serves to clarify much of the early discussion surrounding infinitesimals, enabling a frank discussion of how the real numbers are related to the hyperreal numbers and what the hyperreal numbers are like.

\section{Two specific clarifications of the nature of infinitesimals}

Let us imagine that two further clarifying principles are introduced. First, in order to explain the nature and existence of infinitesimals, our imaginary Leibniz writes:
\begin{quote}
(1) \textit{Every conceivable gap in the numbers is filled by infinitesimals.}
\end{quote}
The gap between $0$ and the positive real numbers, for example, is thus filled with the infinitesimal hyperreal numbers, and similarly there are hyperreal numbers at infinitesimal distance to $\sqrt{2}$ and to $\pi$. This idea, of course, amounts to an incipient form of saturation, expressing how the hyperreal number system transcends the real numbers. Namely, we know now  the hyperreal numbers $\R^*$ of nonstandard analysis are countably saturated, which means that every countably specified gap
 $$x_0\leqslant x_1\leqslant x_2\leqslant \cdots\quad\phantom{< z <}\quad\cdots\leqslant y_2 \leqslant y_1 \leqslant y_0$$
with $x_i<y_i$ is filled by some hyperreal number $z$ strictly between
 $$x_0\leqslant x_1\leqslant x_2\leqslant \cdots\quad < z <\quad\cdots\leqslant y_2 \leqslant y_1 \leqslant y_0$$
Indeed, there will be many such $z$ strictly in the gap, since the gap between the $x_n$ and $z$ itself will also get filled, as will the gap above $z$ and below all the $y_n$.

In our actual history, the actual Leibniz was already inclined toward (1), considering higher orders of infinitesimality. According to~\cite[p.173]{Jesseph1993:Berkeleys-philosophy-of-mathematics}, Berkeley complains:
\begin{quote}\footnotesize
Some mathematicians (notably Leibniz and L'Hopital) hold that there are infinitesimal quantities of all orders and
``assert that there are infinitesimals of infinitesimals of infinitesimals,
without ever coming to an end.''
\end{quote}
The historical Euler 1780 (\cite{Euler1780:On-the-infinity-of-infinities-of-orders}; see also~\cite[p.87--88]{Button+Walsh2018:Philosophy-and-model-theory}) also was busy exploring the vast space of infinite orders of the infinitely large and the infinitely small, discovering the saturation-like manner in which gaps get filled. He observed that for any infinitely large quantity $x$ the number $x^2$ will be infinitely larger still and $x^3$ infinitely larger than that. In contrast, $\sqrt{x}$ will be infinitely smaller than $x$, but still infinite, and the higher roots $\sqrt[3]{x}$, $\sqrt[4]{x}$, similarly get infinitely smaller with each step, while remaining infinite. Nevertheless, Euler demonstrates that we may find orders of infinity that are infinitely smaller than every $\sqrt[n]{x}$, but still infinite. One such number is $\ln x$, and then $(\ln x)^2$ will be infinitely larger than this, but still smaller than every $\sqrt[n]{x}$, and similarly $\sqrt{\ln x}$ smaller again. By taking reciprocals, he finds the same rich phenomenon amongst the orders of infinitesimality relevant for calculus. In this way, he makes a festive party out of filling gaps in the orders of infinity, realizing more and more instances of saturation and thereby supporting principle~(1). And of course this kind of thinking fed into the much later work of Hardy~\cite{Hardy1910:Orders-of-infinity} on the orders of infinity and Hausdorff's related work showing in his context that all countably specified gaps are filled.\goodbreak

Second, in order to justify his calculations with fluxions, the ultimate ratios, and evanescent increments, let us imagine that Newton writes:\nobreak
\begin{quote}
(2) \textit{The two number realms fulfill all the same fundamental mathematical laws.}
\end{quote}\goodbreak

According to the new dictum, the hyperreal numbers would thus fulfill the associativity and distributivity laws, or indeed any law that is true for the real numbers. From our perspective, of course, we can view this statement as an incipient form of the transfer principle, by which the hyperreal field is an elementary extension of the real field $\R\elesub\R^*$, even in expansions of the field structure to include other functions and relations. In particular, $\R^*$ would be a real-closed field, an ordered field in which every positive number has a square root and every odd-degree polynomial has a root.

Jesseph~\cite[p.135]{Jesseph1993:Berkeleys-philosophy-of-mathematics} writes that ``Wallis [1685] took infinitesimal methods to be essentially the same as the method of exhaustion but shorter and more readily applied,'' which can be seen as a conservativity claim about the methods, and Button and Walsh~\cite[\S4.7]{Button+Walsh2018:Philosophy-and-model-theory} describe Leibniz's philosophy of fictionalism about infinitesimals, expressed in his letter to Varignon, which can also be seen as a form of conservativity in that there is nothing really new in them. Newton himself in the \emph{Principia} (1687) expresses conservativity for his methods, stating:
\begin{quote}\footnotesize
These Lemmas are premised to avoid the tediousness of deducing perplexed demonstrations \emph{ad absurdum}, according to the method of the ancient geometers.~\cite[p.102]{Newton:Principia}
\end{quote}
Thus again, there is nothing mathematically new going on. This conservativity attitude fits with the proto-transfer-principle idea of statement (2), which also expresses a kind of conservativity, that there are no new mathematical rules arising in the new hyperreal number realm.

The proposal here with my thought experiment is definitely not that Newton and Leibniz have a full-blown well formulated theory of saturation and transfer, or even of the concept of a field (which only came much later, after Galois), but rather only that they have expressed the primitive idea of two distinct number realms, the real numbers and hyperreal numbers, with vaguely expressed ideas that appear from a contemporary perspective as incipient forms of saturation and the transfer principle. The thought experiment requires only very small initial steps towards the two-realms conception, since further development and rigor would naturally come in time, just as it did in our actual mathematical history.\goodbreak

Nevertheless, the proposal does call for us to imagine that mathematicians might have been a little more modern in their attitude toward number systems, moving beyond the historical understanding of numbers at that time as ratios of geometric magnitudes toward the idea of there being distinct number realms. To be sure, mathematicians have long distinguished at least between natural numbers and other kinds of magnitudes, such as when considering number-theoretic concepts such as even, odd, and prime, and the thought experiment calls for further analogous distinctions concerning the infinitesimals.

Meanwhile, the extent to which the historical Newton had used infinitesimals in the first place is a matter of discussion by historians of mathematics. He had renounced them explicitly, taking himself to mount instead his method of fluxions and the concept of ultimate ratios, which can be seen arguably either as a proto-limit concept or as disguised infinitesimals. In any case, for the purposes of my thought experiment we needn't be troubled by Newton's possible hesitancy towards infinitesimals, since the thought experiment is not about the historical Newton, but about the infinitesimal concept itself, which certainly did exist at the time. For the success of the thought experiment, therefore, if this is an issue we could simply imagine if necessary that it was mainly Leibniz or another Leibniz-like figure who had provided the somewhat fuller explication of the nature of infinitesimals that I am discussing.

\section{The hyperreals are fundamentally coherent}

The main initial point I should like to stress is that we know that these ideas about saturation and the transfer principle for the real and hyperreal numbers are fundamentally coherent and mathematically correct in light of the much later developments of nonstandard analysis, due to Abraham Robinson in the 1960s; they are definitely sufficient for a robust infinitesimal theory of calculus. A glance at Keisler's remarkable infinitesimals-based undergraduate calculus text~\cite{Keisler2000:Elementary-calculus-an-infinitesimal-approach} shows what is possible when beginning even with only very elementary ideas---the entire classical theory can be developed on these notions. To be sure, in our actual mathematical history, the development of calculus proceeded to enormous success on the basis of very primitive infinitesimal ideas, without a fully rigorous foundation, and yet still achieved the key mathematical insights. The thought experiment I propose is that all those developments and insight would still occur, of course, and more, because even a slightly greater initial clarity in the infinitesimal concept would naturally lead to and support fruitful further analysis. In the imaginary history, the development of calculus would be something a little closer to the developments of what we now call nonstandard analysis, perhaps primitive at first, but with increasing sophistication and rigor. Precisely because we know that calculus can be successfully developed this way, the approach would not meet with any fundamental obstacle.

\newpage
Kurt \Godel\ explains his views on nonstandard analysis after Robinson's talk on the subject at the IAS Princeton in 1973:
\begin{quote}\footnotesize
There are good reasons to believe that non-standard analysis, in some version or other, will be the analysis of the future.

One reason is the just mentioned simplification of proofs, since simplification facilitates discovery. Another, even more convincing reason, is the following: Arithmetic starts with the integers and proceeds by successively enlarging the number system by rational and negative numbers, irrational numbers, etc. But the next quite natural step after the reals, namely the introduction of infinitesimals, has simply been omitted.
I think in coming centuries it will be considered a great oddity in the history of mathematics that the first exact theory of infinitesimals was developed 300 years after the invention of the differential calculus.
\cite[p.~311]{Godel:Collected-works-II}
\end{quote}
\Godel\ thus paints the picture that it is our own actual mathematical history that is odd and strange---the history of ideas in my imaginary thought experiment, in contrast, would be the more natural progression. I take this as truly very strong support for the fundamental reasonableness of my thought experiment.

\section{The hyperreal numbers become a familiar mathematical structure}

In this imaginary early history, therefore, the hyperreal numbers would be successful in the foundations of calculus, and as a result they would be taken seriously as a distinct realm of numbers, becoming a core part of the mathematical conceptions underlying the calculus. In the imaginary world, calculus would be founded fully on infinitesimals, without any need for $\forall\varepsilon\,\exists\delta$ limit concepts; the use of infinitesimals would become increasingly sophisticated and rigorous.

The hyperreal numbers would thus enter the Pantheon of number systems at the center of mathematics, the fundamental structures that mathematicians discovered and then used throughout their mathematical work.
$$\N\quad\qquad\Z\quad\qquad\Q\quad\qquad\R\quad\qquad\C\quad\qquad\R^*$$
The hyperreal numbers would thereby find their place in mathematics alongside the other familiar standard mathematical number systems---the natural numbers, the integers, the rational numbers, the real numbers, the complex numbers---and in the world of my thought experiment, there would stand also the hyperreal numbers.

Reflecting on the move from the real numbers $\R$ to the hyperreal numbers $\R^*$, one might be encouraged to seek out saturated versions of all our favored mathematical structures. But actually, for the number systems I have mentioned in the Pantheon above, the hyperreals already provide this. By the transfer principle, we get the hypernatural numbers $\N^*$, the ring of hyperintegers $\Z^*$, the field of hyperrational numbers $\Q^*$, all sitting inside the hyperreal number field $\R^*$, as well as the hypercomplex numbers $\C^*=\R^*\![i]$, consisting of numbers $a+bi$, where $a,b\in\R^*$ are hyperreal. In this way, we may view the move to the hyperreals simply as the move of saturating all our familiar structures.

\section{On the necessity of categoricity for structuralism}

Daniel Isaacson~\cite{Isaacson2011:TheRealityOfMathematicsAndTheCaseOfSetTheory}, taking inspiration from Kreisel, describes the process by which mathematicians come to know their mathematical structures. Namely, as I describe it in~\cite{Hamkins2021:Lectures-on-the-philosophy-of-mathematics}, by informal rigor ``we become familiar with a structure; we find the essential features of that structure; and then we prove that those features axiomatically characterize the structure up to isomorphism. For Isaacson, this is what it means to identify a particular mathematical structure, such as the natural numbers, the integers, the real numbers, or indeed, even the set-theoretic universe.'' Isaacson says,
\begin{quote}\footnotesize
\ldots the reality of mathematics turns ultimately on
the reality of particular structures. The reality of a particular structure,
constituting the subject matter of a branch of mathematics such as
number theory or real analysis, is given by its categorical characterization, i.e. principles which determine this structure to within isomorphism.
\cite[p.~2]{Isaacson2011:TheRealityOfMathematicsAndTheCaseOfSetTheory}
\end{quote}
In this way, the categorical characterizations of our familiar particular structures become the framework of our mathematical reality.

This plays out in our actual mathematical history when Dedekind~\cite{Dedekind1888:What-are-numbers-and-what-should-they-be} proves that the natural number structure $\N$ is uniquely specified up to isomorphism by his theory of the successor operation, leading directly to Peano's elegant development of elementary number theory in that framework. Building upon this, mathematicians provide categorical accounts of the integer ring $\Z$ and the rational field $\Q$. Cantor~\cite{Cantor1895:Beitrage-zur-Begrundung-der-transfiniten-mengenlehre, Cantor1897:Beitrage-zur-Begrundung-der-transfiniten-mengenlehre, Cantor:Contributions-to-the-founding-of-the-theory-of-transfinite-numbers} proves that the rational order $\Q$ is characterized as the unique countable endless dense linear order. Huntington~\cite{Huntington1903:Complete-sets-of-postulates-for-the-theory-of-real-quantities} provides the categorical account of the real field $\R$ as the unique complete ordered field. The complex numbers are characterized as the algebraic closure of $\R$. The categorical characterizations of these core mathematical structures are the central results for the coherence of the mathematical enterprise, enabling us to refer to the various fundamental mathematical structures by their defining characteristics.

I have pointed to the categoricity results as the origin of the philosophy of structuralism in mathematics.
\begin{quote}\footnotesize
Categoricity is central to structuralism because it shows that the essence of our familiar mathematical domains, including $\N$, $\Z$, $\Q$, $\R$, $\C$, and so on, are determined by structural features that we can identify and express. Indeed, how else could we ever pick out a definite mathematical structure, except by identifying a categorical theory that is true in it? Because of categoricity, we need not set up a standard canonical copy of the natural numbers, like the iron rod kept in Paris that defined the standard meter; rather, we can investigate independently whether any given structure exhibits the right structural features by investigating whether it fulfills the categorical characterization. \cite[p.~31]{Hamkins2021:Lectures-on-the-philosophy-of-mathematics}
\end{quote}
In short, having categorical accounts of all our core mathematical structures is necessary for mathematical reference and supports a structuralist mathematical practice.

In set theory, this phenomenon arises after the improvement of Zermelo's flawed initial set theory to Zermelo-Fraenkel set theory, with the addition of the replacement and foundation axioms, an improvement that makes possible Zermelo's famous quasi-categoricity results of~\cite{Zermelo:1930}, showing that the models of second-order set theory $\ZFC_2$ agree with one-another on initial segments. The new $\ZFC_2$ theory thus enjoys a measure of categoricity for the intended set-theoretic universe, leaving open only how high the ordinals will grow. To my way of thinking, categoricity should be a bigger part of the conversation concerning Zermelo's original set theory versus Zermelo-Fraenkel set theory. The models of this theory $\ZFC_2$ are exactly the uncountable Grothendieck-Zermelo universes, now used pervasively in the foundations of category theory, and of course set theorists study them in connection with the inaccessible cardinals. Many of these models admit fully categorical characterizations, and Robin Solberg and I~\cite{HamkinsSolberg:Categorical-large-cardinals} explore the spectrum of fully categorical extensions of $\ZFC_2$, while considering the curious tension between categoricity and reflection principles in the foundations of set theory.\goodbreak

\section{The key imaginary event}

In the world of my historical thought experiment, we shall similarly have categorical characterizations of all the various principal mathematical structures at the end of the 19th century and early 20th century, including the natural numbers, the integers, the real numbers, and the complex numbers.

But what about the hyperreal numbers? In the imaginary history, after all, the hyperreal numbers $\R^*$ have become a core mathematical structure alongside all the others, situated at the very foundations of calculus, present from the start of that subject. Mathematicians would demand an account of the definitive underlying theory of the hyperreal numbers, which would require a categorical characterization like all the others. Naturally, one would expect this characterization to involve the key features already recognized as characteristic of the hyperreal numbers, the saturation ideas and the transfer principle, just as the characterization of the natural numbers involves induction and that of the real numbers involves the least-upper-bound completeness principle.

Thus we come to the key imaginary event of the thought experiment. Namely, let us imagine that in the early 20th century, a Zermelo-like figure formulates a sufficient theory---the theory $\ZFC+\CH$ suffices---able to prove a categorical characterization of the hyperreal number field $\R^*$, similar to how the actual Zermelo introduced his set theory as an explanation of his proof of the well-order theorem.

The theory $\ZFC+\CH$ is indeed able to provide the desired characterization of the hyperreal field $\R^*$ as follows, which shows under~\CH\ how $\R^*$ is characterized by refined versions of the two ideas we had attributed to Newton and Leibniz in the thought experiment.
\newtheorem*{hrct}{Hyperreal categoricity theorem}
\begin{hrct}
 Assume $\ZFC+\CH$. Then there is up to isomorphism a unique smallest countably saturated real-closed field.
\end{hrct}
This theorem is by now a standard result in elementary model theory, proved using a back-and-forth argument in the style of Cantor's famous argument about the rational order, except that here the back-and-forth construction proceeds transfinitely through $\omega_1$ many steps rather than just countably many (see Erd\H{o}s, Gillman, and Henriksen~\cite{ErdosGillmanHenriksen1955:An-isomorphism-theorem-for-real-closed-fields}). A model is \emph{countably saturated}, that is, $\aleph_1$-saturated, if it realizes every finitely satisfiable type with countably many parameters, but in the theorem we actually need to require only that the order is saturated---all countably described gaps should be filled. The general fact in play here is that any two saturated models of the same complete theory and size are isomorphic, and that would be the situation for our hyperreal fields under the hypotheses of the categoricity theorem. The smallest possible size in question would be the continuum, since countably saturated real-closed fields must have size at least continuum, and there are such fields of size continuum.

Indeed, we have several various constructions of a countably saturated real-closed field of size continuum. Namely, (1) we can construct the ultrapowers $\R^\N/\mu$ of the real field by any nonprincipal ultrafilter $\mu$ on $\N$, and this is always a countably saturated real-closed field of size continuum; (2) we can proceed via Hahn series, a generalized kind of power series, to construct a countably saturated real-closed field of size continuum; (3) we can undertake a general model-theoretic construction, successively realizing types in a transfinite elementary chain, to produce a countably saturated model of any given consistent theory, including the theory of real-closed fields; (4) perhaps exemplifying this construction in an attractive, concrete general manner, we can undertake the Conway construction of the surreal field through all countable ordinal birthdays, filling all possible gaps that arise---the result is $\No(\omega_1)$, a countably saturated real-closed field of size continuum.

The hyperreal categoricity theorem shows under \ZFC+\CH\ that all these various constructions give rise to exactly the same hyperreal field, which we may thereby regard as the canonical structure of the hyperreal field $\R^*$. The situation for the hyperreals under \CH\ is thus rather like that of the real field in \ZFC, for which we also have a variety of constructions, proceeding with Dedekind-cuts in the rationals or equivalence classes of Cauchy sequences and so forth. All the various presentations of complete ordered fields are provably isomorphic in \ZFC, and this categoricity enables a structuralist treatment of the real field as the unique complete ordered field. In a sense \ZFC\ is aimed at providing the satisfactory theory of the real numbers, for the real categoricity theorem is not provable in weaker systems, such as constructive mathematics, where mathematicians must treat the Dedekind reals as a distinct conception from the Cauchy reals.

Similarly, in \ZFC+\CH, all the various constructions of the hyperreal field give rise to the same underlying canonical structure, thereby enabling a structuralist account of infinitesimals---the hyperreals are the unique smallest countably saturated real-closed field.

\section{\CH\ is required}\label{Section.CH-required}

I should like to call attention to a key feature of the thought experiment, namely, the mathematical fact that the~\CH\ is required. We can provide the categorical characterization of the hyperreal numbers in \ZFC+\CH\ as stated in the hyperreal categoricity theorem, but this is not possible in $\ZFC$ alone. Judith Roitman~\cite{Roitman1982:Non-isomorphic-hyper-real-fields} showed that it is relatively consistent with $\ZFC+\neg\CH$ that there are multiple non-isomorphic hyperreal fields arising as ultrapowers $\R^\omega/\mu$, and these are always countably saturated real-closed fields of size continuum. Alan Dow~\cite{Dow1984:On-ultrapowers-of-Boolean-algebras} showed that whenever~\CH\ fails, then indeed there are multiple non-isomorphic ultrapowers $\R^\omega/\mu$, non-isomorphic even merely in their order structure. Thus,~\CH\ is outright equivalent to the hyperreal categoricity assertion that there is a unique smallest countably saturated real-closed field (see also \cite{Esterle1977:Solution-to-problem-of-Erdos-Gillman-Henriksen}).

These results show that the phrase ``the hyperreal numbers'' is not generally meaningful in \ZFC, because in $\ZFC$ there is not necessarily just one mathematical structure fitting the description. But in \ZFC+\CH, there is. With the continuum hypothesis, we can specify the hyperreal numbers up to isomorphism as a canonical structure, the unique smallest countably saturated real-closed field.

\section{How~\CH\ gets on the list}

So this is how~\CH\ gets on the list of fundamental axioms. The thought experiment, at bottom, is that the hyperreal field $\R^*$ is long a core mathematical idea, pre-rigorous at first, but then with increasing rigor and sophistication. To give a foundational account of the hyperreal number system and thus of the theory of infinitesimals, the Zermelo-like figure provides an existence proof and categorical characterization, introducing the fundamental axioms of $\ZFC+\CH$ in order to do so. We know that this is possible and, furthermore, that~\CH\ cannot be omitted. So~\CH\ gets onto the list of fundamental axioms, being necessary to establish the basic coherence or even (in the Isaacson sense) the reality of the hyperreal numbers and thus indispensable for the foundations of calculus.

\section{Extrinsic and intrinsic support for~\CH}

The developments would provide enormous extrinsic support for~\CH, similar to the extrinsic justification $\ZFC$ currently enjoys in light of its robust foundational account of the real numbers $\R$. The theory $\ZFC+\CH$ would be seen as similarly successful regarding the theory of the hyperreal numbers.

In the imaginary history, I would find it quite likely that after the~\CH\ had found its extrinsic justification in this way as a mathematical necessity, then intrinsic justifications would also begin to find their appeal, similar to how the axiom of choice is often viewed as extrinsically justified by its widespread use and important central consequences, but set theorists also point to its intrinsic justification under the concept of arbitrary set existence. In the case of the continuum hypothesis, the intrinsic justification I imagine is that~\CH\ asserts that the two methods of achieving uncountability agree, that is, the process of going to the next higher cardinal gives the same result as taking the power set; in short, $\aleph_1=\beth_1$. This can be seen as a unifying, explanatory principle of the uncountable, and therefore an intrinsic justification for~\CH.\goodbreak

We might also reflect on the fact that by basic human rationalizing nature, one naturally finds it easier to be convinced by arguments for the intrinsic truth of an axiom, once one has already been convinced of the axiom's extrinsic necessity.

\section{A generalized thought experiment and the generalized continuum hypothesis}\label{Section.GCH}

The hyperreal categoricity theorem is that under~\CH, the hyperreal field is the unique smallest countably saturated real-closed field. Perhaps a critic objects that the smallest-size requirement is ad hoc---wouldn't it be more natural to relax this and consider hyperreal fields of other sizes? And isn't it unnatural to require only countable saturation, rather than full saturation?

Let me address this initially by defending countable saturation as a rich, natural notion. Countable saturation is both easy to express and understand, and suffices for the robust existence of infinitesimals in a vast hierarchy of orders. It leads to a rich, successful theory, without the need for higher levels of saturation. Furthermore, the construction of a countably saturated real-closed field, as with the surreal construction through the countable birthdays, is both clear and natural. So there is very little lacking in our conception of the hyperreals as a countably saturated real-closed field. And Hausdorff's 1909 proof that $\N^\N/\text{Fin}$ is countably saturated but has an unfilled $(\omega_1,\omega_1)$ gap would tend to encourage a greater focus specifically on countable saturation.

Meanwhile, to be sure, full saturation does imply categoricity in any given cardinality in which it occurs by the back-and-forth construction. And furthermore, the existence of fully saturated models does not require~\CH. The existence of saturated real-closed field of size continuum is equivalent merely to $\continuum^{<\continuum}=\continuum$, which can occur say, even if $\continuum=2^{\aleph_0}=2^{\aleph_1}=\aleph_2$, and in many other kinds of cases. There will be a fully saturated real-closed field of uncountable size $\kappa$ if and only if $\kappa^{<\kappa}=\kappa$, which is independent of the \CH\ and \GCH, and it can hold consistently by forcing with any particular regular uncountable cardinal, although it holds necessarily of any inaccessible cardinal. And yet, even in these cases, where \CH\ fails and there is a saturated field of size continuum or more, there will also be numerous non-isomorphic countably saturated such fields, which may detract from the sense of uniqueness for the hyperreal conception.

Nevertheless, let me take on board both of the critical objections at once with a generalized thought experiment. Let me imagine that the Zermelo-like figure adopts an expansive attitude toward the hyperreals, proving instead the following generalized hyperreal categoricity theorem, on the basis of the \emph{generalized continuum hypothesis} (\GCH), which asserts that $2^\kappa=\kappa^+$ for every infinite cardinal~$\kappa$.

\newtheorem*{ghrct}{Generalized hyperreal categoricity theorem}
\begin{ghrct}
Assume $\ZFC+\GCH$. Then there are up to isomorphism unique saturated real-closed fields in every uncountable regular cardinality.
\end{ghrct}

The theorem can be proved by observing that saturated models of the same size and theory are unique up to isomorphism by the back-and-forth construction, as we have mentioned, and you get existence of saturated models from the \GCH\ by realizing types in a transfinite elementary tower. In fact the converse of this theorem also is true---that is, the \GCH\ is required for the generalized categoricity here---since $\kappa^{<\kappa}=\kappa$ is necessary for there to be a $\kappa$-saturated model of size $\kappa$, and this implies the~\GCH\ if it holds for every uncountable regular cardinal. Thus, the generalized hyperreal categoricity principle is itself equivalent over \ZFC\ to the~\GCH.

On the basis of the \GCH\ we thereby find ourselves blessed with canonical hyperreal fields in every desired size, a transfinite tower of orders of infinitesimality continuing to all higher cardinals. This would be a natural continuation of the saturation ideas originating with Leibniz and continuing with Euler's exploration of the diverse orders of infinity, and then with Hardy and into contemporary times. And since the \GCH\ is required, this is how the \GCH\ could also get on the list of fundamental axioms. These higher uncountable hyperreal fields would be seen as converging in a vast elementary chain ultimately to the surreal numbers, another core number concept with a proper-class categoricity characterization.

\section{Foundationalism and a priori knowledge}

Williamson~\cite[\S III]{Williamson2016:Absolute-provability-and-safe-knowledge-of-axioms} speculates on how it would be that other beings, with physical and mathematical powers similar to humans, might settle the continuum hypothesis.
\begin{quote}\footnotesize
One normal mathematical process, even if a comparatively uncommon one, is adopting a new axiom. If set theorists finally resolve~\CH, that is how they will do it. Of course, just arbitrarily assigning some formula the status of an axiom does not count as a normal mathematical process, because doing so fails to make the formula part of mathematical knowledge. In particular, we cannot resolve~\CH\ simply by tossing a coin and adding~\CH\ as an axiom to ZFC if it comes up heads, ${\sim}\CH$ if it comes up tails. We want to know whether~\CH\ holds, not merely to have a true or false belief one way or the other (even if we could get ourselves to believe the new axiom). Thus the question arises: when does acceptance of an axiom constitute mathematical knowledge?
\end{quote}
My thought experiment engages with Williamson's challenge by describing the richer context and process that would lead mathematicians to the~\CH. In the alternative world I describe, mathematicians have gained increasing familiarity over the centuries with the hyperreal number system, embedding it at the center of classical mathematical developments, especially in the calculus, and thus they have gained increasing confidence in their use of the hyperreals as a particular, familiar mathematical structure. Seeking a more thorough underlying mathematical explanation of it, including a categorical account, as we do with all our particular mathematical structures, they would find it in a theory including~\CH, which can prove hyperreal categoricity, and we know that~\CH\ is necessary for this. In this way, my thought-experiment mathematicians are led to the~\CH\ by undertaking what Williamson calls the ``normal mathematical process.'' The axiom they accept in effect is hyperreal categoricity, since it is necessary for the coherence and reality of their mathematical practice, and this principle is equivalent to the~\CH.

In another thought experiment, Berry~\cite{Berry2013:Default-reasonableness-and-the-mathoids} introduces the mathoids, creatures who look upon Fermat's last theorem as intuitively and immediately obvious, a foundational belief, without feeling any need to justify this stance with an argument that we would find convincing; similar imaginary creatures are discussed in~\cite[\S III]{Williamson2016:Absolute-provability-and-safe-knowledge-of-axioms}. This is part of her investigation into foundationalism about a priori mathematical knowledge, in which she highlights a problematic conclusion, namely, ``that nothing in the current literature lets us draw a principled distinction between what these creatures are doing and paradigmatic cases of good a priori reasoning.''

I take myself to sidestep that particular debate in my thought experiment, precisely because I took pains to describe mathematical characters that we \emph{would} find convincing. My imaginary Newton and Leibniz and the imaginary Zermelo-like figure seem perfectly reasonable given what we know about the underlying mathematics---far more reasonable than the mathoids, who strike us (and this is important for Berry's point) as unreasonable in their mathematical beliefs. My thought experiment, after all, is simply to reorder certain mathematical discoveries that we know are correct, in such a way that the hyperreals would naturally become far more central in classical mathematics than they are in our actual history, with the consequence that the need for a categorical account of them would become more urgent. This would lead us, I have argued, to a very different attitude towards the foundational theories able to provide such an account, and these theories would have to include~\CH.\goodbreak

\section{Imaginary later history}

Let us continue the thought experiment by considering how later mathematical developments would be received in the world of the imaginary history I have proposed. \Godel\ proves that $\ZFC+\CH$ is true in the constructible universe $L$, which would of course be very welcome confirmation of the main theory. This could even lend some extrinsic support for $V=L$, more so than it enjoys currently. Indeed, the fact that \GCH\ holds in $L$ would provide further support for $V=L$ in light of the generalized hyperreal categoricity result in section~\ref{Section.GCH}. Solovay's theorem~\cite{Solovay2006.FOM:V=L-implies-complete-second-order-theories-are-categorical}, showing that under $V=L$ every finitely axiomatizable complete second-order theory is categorical, might be welcomed as a corresponding fundamental principle in the same light, rather than a curiosity about $L$ as currently. (Solovay's theorem has been generalized to $L[\mu]$ and the large cardinal context by~\cite{SaarinenVaananenWoodin:On-the-categoricity-of-complete-second-order-theories}, who also show under the axiom of projective determinacy, a unifying consequence of large cardinals, that every finitely axiomatizable complete second-order theory with a countable model is categorical.)

A commitment to a robust theory of the hyperreals presumes in certain ways a commitment to the axiom of choice or at least fragments of it. One needs the prime ideal theorem to have suitable ultrafilters $\mu$ with which to form the ultrapower $\R^\N/\mu$, and one needs countable choice in order to know that indeed these satisfy the transfer principle and are countably saturated. Conversely, the existence of a hyperreal field with infinitesimals and the transfer principle outright implies the existence of nonprincipal ultrafilters on the natural numbers, since for any fixed infinite number $N\in\N^*$ we can define $X\in\mu\iff X\of\N\text{ and }N\in X^*$, in effect defining that $X$ is $\mu$-large when it expresses a property that $N$ exhibits. Meanwhile, the existence of ultrafilters implies the existence of non-Lebesgue measurable sets in the real numbers, an often-mentioned consequence of the axiom of choice. In this way, my thought-experiment inhabitants take the \emph{modus tollens} to Alan Connes's criticism of nonstandard analysis (see \cite{GoldsternSkandalis2007:Interview-with-A-Connes} and my related discussion \cite[p.~81]{Hamkins2021:Lectures-on-the-philosophy-of-mathematics}). In this way, a commitment to the hyperreal numbers carries a small accompanying commitment to the axiom of choice.

Meanwhile, the discovery via forcing that without~\CH\ there can be multiple non-isomorphic hyperreal fields, as mentioned in section~\ref{Section.CH-required}, would be seen as chaotic and bizarre, perhaps a little like current attitudes about models of \ZF\ with strange failures of the axiom of choice. For example, it is known to be relatively consistent with \ZF\ without the axiom of choice that the rational field $\Q$ can have multiple distinct non-isomorphic algebraic closures, a countable one as well as an uncountable one~\cite{Lauchli1962/63:Auswahlaxiom-in-der-algebra,Hodges1976:Lauchlis-algebraic-closure-of-Q}. Mathematicians often find this situation very strange, and many regard this model as getting something fundamentally wrong about the algebraic numbers. I believe similarly that the mathematicians in the imaginary world would find it very odd to have multiple non-isomorphic hyperreal fields, and this would reinforce the view that $\ZFC+\CH$ is the right theory.

The method of forcing would be received differently in the imaginary world than in our own world---it would be perceived as a little less successful. For us, one of the attractive central features of forcing is that it necessarily preserves \ZFC. Every forcing extension of a model of $\ZFC$ is another model of \ZFC. But in the imaginary world, the main standard theory is \ZFC+\CH, and the corresponding feature is not true of this theory, since forcing can destroy~\CH. Indeed, this was Cohen's main initial application, to produce a model of $\ZFC+\neg\CH$. Perhaps the forcing method would be viewed in a similar way to how some people currently view the symmetric model constructions, which preserve \ZF, but not necessarily the axiom of choice. We often use the symmetric model construction to build strange badly behaved models of \ZF, enabling us to see how things can go awry when one doesn't have the axiom of choice. Similarly, in the imaginary world I propose, forcing could be used to build ``strange'' models of $\ZFC$ with $\neg\CH$, showing how things can go awry when one doesn't have the continuum hypothesis.

In particular, my imaginary mathematicians would have a very different attitude concerning the dream-solution situation I described in the introduction, since in the thought experiment, they have adopted~\CH\ as fundamental before having had the experience via forcing of the $\neg\CH$ worlds, which they would view as strange. My thought experiment shows in effect how it could have been that the dream solution was achieved, even though it is no longer possible for us to achieve it. For my imaginary mathematicians,~\CH\ is a consequence of categoricity for the hyperreal numbers, and that is for them an instance of the dream solution.

In the imaginary world, of course eventually the $\forall\varepsilon\,\exists\delta$ accounts of limits and continuity would be discovered, but these would be seen as complicated and unnecessary abstractions, in light of the straightforward use of infinitesimals. The situation would be an inversion of current attitudes towards nonstandard analysis and ultrapowers.

A challenging counterpoint would eventually flow from the discovery by Esterle \cite{Esterle1977:Solution-to-problem-of-Erdos-Gillman-Henriksen} (see also \cite[theorem~17]{Ehrlich2012:Absolute-arithmetic-continuum-and-the-unification-of-all-numbers-great-and-small}) that in \ZFC\ there is an initial countably saturated real-closed field, that is, one that embeds isomorphically into all other such fields, without requiring \CH. By stratifying the field as a union of fields with all cuts having countable cofinality, one can mount a subtle back-and-forth argument along that hierarchy to show that this field is unique, again without \CH. Indeed, this initial field is simply the field $\No(\omega_1)$ of the surreal numbers born at a countable ordinal birthday. This is a highly canonical mathematical structure, a countably saturated real-closed field that admits a very natural construction and furthermore embeds isomorphically into all such fields---it is therefore ``smallest'' in a stronger way than mere cardinality as in the hyperreal categoricity theorem, and furthermore its categorical characterization does not require \CH. Would this undermine the status of the hyperreal categoricity theorem in my imaginary world? Perhaps people would want to replace the previous understanding of what the hyperreal numbers are with a new categorical account---they are the unique initial countably saturated real-closed field.

I would have several responses to this. First, the situation when \CH\ fails, as I have mentioned, would remain chaotic, with numerous non-isomorphic minimal-size countably saturated real-closed fields, arising as ultrapowers $\R^\N/\mu$ and by other constructions, including $\No(\omega_1)$. Second, while the initiality of $\No(\omega_1)$ is a very natural property, I totally agree, meanwhile there would be other hyperreal field candidates with other very natural properties. For example, if $\continuum^{<\continuum}=\continuum$ there would also be a fully saturated real-closed field of size continuum, which would definitely be different than $\No(\omega_1)$ when \CH\ fails, and yet this alternative would also be a natural, canonical structure. Which would be the real hyperreals? It recalls the situation of the real numbers in constructive mathematics, where one must distinguish between the Dedekind reals and the Cauchy reals and so forth, and these are not constructively isomorphic. Which are the real reals in constructive mathematics? Similarly, which are the real hyperreals in $\ZFC+\neg\CH$? Under \CH, all the various candidates are isomorphic, which largely dissolves the issue, and in this sense, the categoricity situation is simply much better with \CH\ and better still with \GCH. But meanwhile, perhaps in my imaginary world there would be a community of mathematicians proving theorems about the various kinds of hyperreal fields in mere \ZFC, in the same way that we currently have a community of mathematicians proving theorems about the various real fields in constructive mathematics.

\section{Hyperreal hesitancy in the actual world}

To my way of thinking (see~\cite[p.~78-79,~Questions~2.18,~2.19,~2.20]{Hamkins2021:Lectures-on-the-philosophy-of-mathematics}), the lack of a categoricity result for the hyperreal numbers in $\ZFC$ is a principal part of the explanation for the hesitancy amongst many mathematicians to take up nonstandard analysis. Mathematicians are naturally loathe to mount a fundamental mathematical theory like calculus with an underspecified mathematical structure at its core. We don't want to found calculus on an unknown hyperreal structure, if there are multiple non-isomorphic hyperreal structures to choose from. Which one do we choose? What if the resulting theory and features would depend on the particular choice? And how are we to refer anyway to a particular one of the hyperreal fields in the first place, if we lack a categorical description that picks it out from the alternatives? The need for categorical characterizations of all our fundamental mathematical structures seems necessary for a coherent structuralist mathematical practice of referring to them.

This explanation of the observed hyperreal hesitancy amounts to the contrapositive of the claim that we require categorical accounts for all our core mathematical structures. In our actual mathematical history, $\ZFC$ is already established as the core foundational theory before hyperreals are discovered, and so the lack of categoricity for the hyperreals means that we cannot accept them as a core structure.

One might inquire whether the proposed argument for~\CH\ could be convincing in our actual world, now that we know about the hyperreal number systems. One problem is that in our actual history, the hyperreal numbers never became a core mathematical structure in the first place, and so there is perhaps no pressing need to provide a categorical account of them. In actual history, we made calculus rigorous with the $\forall\varepsilon\,\exists\delta$ limit formalism of Bolzano and Weierstrass, and the hyperreal number systems were a later discovery, a logical curiosity, superfluous for calculus. This is why the hyperreal-categoricity argument for~\CH\ is more compelling in the imaginary thought-experiment world, where hyperreals were a core structure from the start, and less so in the actual world, where they were not.

\section{Historical contingency}

Penelope Maddy~\cite{Maddy1988:BelievingTheAxiomsI} argued that there is a certain historical contingency to the $\ZFC$ axioms of set theory.
\begin{quote}\footnotesize
The fact that these few [\ZFC] axioms are
commonly enshrined in the opening pages of mathematics texts should be viewed
as an historical accident, not a sign of their privileged epistemological or metaphysical status.
\cite{Maddy1988:BelievingTheAxiomsI}
\end{quote}
She was concerned mainly with the large cardinal extensions of $\ZFC$ and the accompanying determinacy principles, seeking reasons to justify their incorporation into our basic conception of set theory.

I have argued in this paper similarly for the historical contingency of $\ZFC$, describing how it could have been that we view the continuum hypothesis itself as a basic principle, necessary for the success of mathematics. While this kind of contingency for \ZFC\ is a consequence of my argument, my main conclusion is not about contingency as such, but rather specifically that we might easily have had very different views about the continuum hypothesis.

Nevertheless, perhaps there would be other historical thought experiments by which we might have come to view $\neg\CH$ as fundamental, although the current examples I know of strike me as less compelling than the example I have described in this paper. Solovay (in personal discussions) has defended a vision of the real continuum that it should be a real-valued measurable cardinal, and this would necessarily involve an outsized failure of the~\CH, by which the continuum is extremely large, but also there would be a fully saturated hyperreal field of size continuum. Moore~\cite{Moore2010:The-PFA} has defended the vision of set theory under the forcing axiom \PFA, and others with \MM\ and $\MM^+$, all of which imply the continuum is $\aleph_2$ and that there is a fully saturated hyperreal field of size continuum. Moore says, ``Forcing axioms have proved very effective in classifying and developing the theory of objects of an uncountable or non separable nature,'' and it would seem possible to expand this with a similar kind of thought experiment by which \PFA\ was essential for the resolution of some fundamental mathematical commitment. Similarly, I have described in~\cite{Hamkins2003:MaximalityPrinciple, HamkinsWoodin2005:NMPccc} forcing axioms such as the c.c.c.~maximality principle $\MP_{ccc}$ and a generalization to the necessary c.c.c.~maximality principle $\necessary\MP_{ccc}$, which imply that the continuum is larger than any cardinal that we can describe in any way that would be absolute to c.c.c. forcing extensions. I can imagine similar thought experiments by which these principles could be taken as fundamental.

If indeed such thought experiments are possible, or indeed only on the basis of the thought experiment of this paper in comparison with the standard \ZFC-only foundations, I am inclined to take them all as an argument against the view that mathematical foundations have a necessary nature. Rather, for various sound reasons, mathematicians might have come to various different and perhaps incompatible conclusions about what they will take to be the central mathematical principles around which they intend to organize their mathematical investigations. For this reason, I regard my thought experiment here as supporting pluralism in the foundations of mathematics, by showing how it could naturally have been that we take a different theory as fundamental than we do currently.

Nevertheless, I recognize that for set theorists taking the universe view, by which there is a unique determinate set-theoretic reality, to have an argument that \CH\ could have been a fundamental truth is in effect to have an argument that it is in fact a fundamental truth. And for this, I shall simply place my thought-experiment argument alongside the other \CH\ arguments on offer. Namely, taken this way, as a proposal to those with the universe view, my argument is that \CH\ is true because it provides the categorical theory of the infinitesimal numbers and indeed it is necessary for this.


\section{Conclusion}

I have described how we could have come to have a very different perspective on the continuum hypothesis. It could easily have been that the early theory of calculus had been a little more clear about infinitesimals, positing that they inhabit a distinct further realm, a system of numbers we might call the hyperreal numbers. This hyperreal number system would thereby have come to be embedded as a necessary component at the core of calculus.
Clarifying the relation between the real and hyperreal numbers, early incipient forms of saturation and the transfer principle would have been put forth, vaguely at first, but then with increasing sophistication. We know from nonstandard analysis that these ideas are capable of serving as foundational in the development of infinitesimals-based calculus, and so the resulting theory would have been robust and successful. With the rise of rigor at the end of the 19th and early 20th centuries, when mathematicians were providing categorical accounts of all the familiar mathematical structures, a Zermelo-like figure would have provided such a characterization of the hyperreal numbers. The theory would be something like \ZFC+\CH, which we know suffices, and also we know that~\CH\ cannot be omitted for this. This would have provided enormous extrinsic justification for the continuum hypothesis, for this axiom would be seen as necessary for making sense of one of the core number systems underlying calculus. Thus, we would view~\CH\ as a fundamental principle necessary for mathematics and indispensable in the foundations of calculus.

To be sure, I am not arguing that~\CH\ already is or should be considered this way as fundamental, but rather only that it could have been. Thus, I claim, we must face a certain degree of contingency in our fundamental theories. What we consider to be bedrock foundational principles could have been different---we could have seen the continuum hypothesis as fundamental.

\AtNextBibliography{\small}
\printbibliography

\end{document}